\documentclass[12pt]{article}

\usepackage{amsmath}

\def \R{I\!\!R}

\def \Z{Z\!\!\!Z}
\def \esssup{\rm{esssup}}

\def\INTER{\mathop{\rm {\cap}}\limits}
\def\ben{\begin{enumerate}}
\def\een{\end{enumerate}}
\def\bit{\begin{itemize}}
\def\eit{\end{itemize}}
\newtheorem{definition}{Definition}
\numberwithin{definition}{section}
\newtheorem{theorem}{Theorem}

\newtheorem{corollary}[theorem]{Corollary}
\newtheorem{proposition}[theorem]{Proposition}
\newtheorem{remark}[theorem]{Remark}
\newtheorem{notation}[theorem]{Notation}
\numberwithin{theorem}{section}

\begin{document}

\title{Bid-Ask  Dynamic Pricing  in Financial  Markets \\
with Transaction Costs and Liquidity Risk}

\author{Jocelyne Bion-Nadal \\CMAP \\{\em UMR 7641 CNRS Ecole Polytechnique}\\ Ecole Polytechnique F-91128 Palaiseau cedex    \\bionnada@cmapx.polytechnique.fr}
\date{ }

\maketitle

\section*{Abstract}
 We introduce, in continuous time, an axiomatic approach to assign to
 any financial position a 
dynamic ask (resp. bid) price process. Taking into account both transaction costs and liquidity risk this leads to the convexity (resp. concavity) of the ask (resp. bid) price. Time consistency is a crucial property for dynamic pricing. Generalizing the result of Jouini and Kallal, we prove that the No Free Lunch condition for a time consistent dynamic pricing procedure (TCPP) is equivalent to the existence of an equivalent  probability measure $R$ that  transforms a process between the  bid process and the ask process of any financial instrument into a martingale. Furthermore  we prove   that the ask price process associated with any financial instrument is then a  $R$-supermartingale process which has a cadlag modification. 
Finally we show that time consistent  dynamic pricing    allows both to extend the dynamics of some reference assets and  to be consistent with any observed bid ask spreads that one wants to take into account.  It then provides new bounds reducing the bid ask spreads for the other financial instruments.
\section*{Introduction} 
\label{sec:intro}  
 
 The theory of asset pricing, and the fundamental theorem, have been formalized by Harrison and Kreps \cite{HK} Harrison and Pliska \cite{HP},  Kreps \cite{KR},  Dalang et al \cite{DMW}, according to no arbitrage principle. In the classical setting, the   market is assumed to be frictionless and then a no arbitrage dynamic price process is a martingale under a probability measure equivalent to the reference probability measure.\\
 However real financial markets are not frictionless and therefore an important literature on pricing under transaction costs  has been developed. Pioneer works  in financial markets with transaction costs are those  of Jouiny  and Kallal  \cite{JK} and  Cvitanic and Karatzas \cite{CK}.
 A matrix formalism introduced by Kabanov \cite{Ka}, for  pricing under proportional transaction costs has been developed in many papers, as in Kabanov and Stricker \cite{KaS}, in Schachermayer \cite{Sc}... 
In  these papers the bid ask spreads are explained by transaction costs. 
 On the contrary the approach developed by Jouini and Kallal \cite{JK}  and extended in Jouini \cite{J2} is an axiomatic approach in continuous time assaigning to financial assets a dynamic ask price process (resp. a dynamic bid price process) in such a way that the ask price procedure is sublinear. Jouini and Kallal \cite{JK}  have proved  that the absence of arbitrage in then equivalent to the existence of an equivalent probability measure which transforms some process between the bid and the ask price process into a martingale. The bid-ask spread in that model can be interpreted as transaction costs or as the result of buying and selling orders.\\
In recent years, a pricing theory has also been developed, taking inspiration from the theory of risk measures. 
This has first been done in a static setting, as in Carr, Geman and Madan \cite{CGM}, F\"ollmer and Schied \cite{FS02}, Staum  \cite{S}, Barrieu and El Karoui \cite{BK} and Bion-Nadal \cite{BN02}. The point of view of pricing via risk measures has also been considered in a dynamic way using Backward Stochastic Differential Equations as in El Karoui and Quenez \cite{EKQ},  El Karoui and Quenez \cite{EKQ2}, El Karoui, Peng and Quenez \cite{EKPQ}, Peng \cite{Peng}, Kl\"oppel and Schweizer \cite{KS} and  Peng \cite{Peng06}. The axiomatization of dynamic pricing procedures goes back to Peng \cite{Peng} where it is called ``g-expectations''. The close notion of monetary utility functions  for processes in a discrete time setting has been studied by  Cheridito, Delbaen and Kupper \cite {CDK}.
A close definition, under the name of Monetary Concave Utility Functional, can be found in Kl\"oppel and Schweizer \cite{KS}. Pricing via dynamic risk measures has also been studied by Jobert and Rogers \cite{JR} in discrete time over a finite probability space. For papers on dynamic risk measures we also refer to Delbaen \cite{D} for the study of coherent (or homogeneous) dynamic risk measures, to Detlefsen and Scandolo \cite{DS} and F\"ollmer and Penner \cite{FP} for dynamic risk measures in a discrete time setting and to Bion-Nadal \cite{BN03} and \cite{BN04} for dynamic risk measures in a continuous time setting.\\
Another approach called good deal bounds pricing has been introduced simultaneously by  Cochrane and Saa-Requejo \cite{CSR} and Bernardo and Ledoit \cite{BL}. The idea is to exclude some payoffs which are too much attractive, called good deals, in  order to restrict the bid ask spread. This corresponds to restrict the set of equivalent martingale measures that one can use to price.  Therefore it is related to the theory of pricing via risk measures. The link, in static case, has been first established by Jaschke and K\"uchler \cite{JaK}. For papers along these lines we refer to Cherny \cite{Ch},
Staum \cite{S},  Bj\"ork and Slinko \cite{BS}, and Kl\"oppel and Schweizer \cite{KS2}.\\
In real financial markets, market partipants submit offers to buy or sell a fixed number of shares. Therefore one can observe the resulting book of orders for any traded security.   The book of order of a financial product $X$ can be interpreted as bid and ask prices associated with the random variables $nX$ where $n$ is a positive integer. The observed ask price associated with $nX$ for $n$ large enough is greater that $n$ times the observed ask price of $X$. This corresponds to the liquidity risk.  Cetin et al \cite{CJP} and \cite{CJPW} or  Astic and Touzi \cite{AT} have developed a pricing  approach taking into account the liquidity risk. 
\\[2mm]
 In this paper we develop an approach which is close from that of Jouini and Kallal \cite{JK}, based on the construction of both bid price process and ask price process in a model free setting. However the question we address here is to define a global dynamic pricing procedure in order to assign  to any  financial instrument (asset, option, basket, portfolio...) a dynamic bid (resp. ask)  price process, taking into account both transaction costs and liquidity risk. 
Assuming that there is at least one liquid asset which is positive, one can take it as num\'eraire. Taking furthermore into account the lack of liquidity, the usefulness of diversification, the existence of transaction costs, it is natural to  impose that the ask price is convex and translation invariant. 
Another property that we require for a dynamic pricing procedure is the time consistency which means that  the price at one instant of time  of a financial instrument can be computed indifferently directly or in two steps using an intermediate instant of time.  We can then take advantage of the theory of time consistent dynamic risk measures, in particular the existence of a dual representation \cite{DS} and the characterization of the time consistency in terms of a cocycle condition \cite{BN03} and \cite{BN04}.\\ The main differences with the paper of Jouini and Kallal are the following: here we don't assume that the (ask) price function is sublinear as in \cite{JK} but we allow it to be convex, in order to take into account the liquidity risk; our construction is a global one in the sense that it assigns to any financial instrument represented by a bounded ${\cal F}_{\tau}$-measurable functional a  dynamic ask (resp bid) price process $(\Pi_{\sigma,\tau}(X))_{\sigma \leq \tau}$ (resp $(-\Pi_{\sigma,\tau}(-X))_{\sigma \leq \tau}$); and we do not assume a priori that the ask (resp bid) price process associated to a given financial instrument has right continuous paths. As in \cite{JK} we define an intrinsic  notion of No Free Lunch for a dynamic pricing procedure. Our main result stated in section \ref{sec:No Free Lunch} is the extension of the theorem proved in \cite{JK} characterizing the No Free Lunch condition.  We prove that the No Free Lunch property is equivalent to the existence of (at least) an equivalent probability measure $R$ which transforms a process between the  bid process and the ask process of any financial instrument into a martingale. This probability measure $R$ is independent of the financial instrument considered. 
Furthermore from   our previous paper on  time consistent dynamic risk measures \cite{BN04} we derive the fact that for any No Free Lunch pricing procedure, for any financial instrument the dynamic bid (resp ask ) price process has a cadlag modification which is a $R$ submartingale (resp supermartingale) process.\\ 
 The last section of the paper deals then with the usual setting for financial markets. The dynamics of some  reference assets $(S^k)_{0 \leq k \leq p}$ are given and we consider time consistent dynamic pricing procedures extending the dynamics of these reference assets. We prove that the upperbound for the ask prices obtained using all such procedures is the surreplication price, generalizing thus to the general convex case the result known for sublinear pricing procedures.  We consider then two  approaches already studied in the literature  to reduce the bid ask spreads which can be formulated in our setting: No good deal pricing and pricing under convex constraints. We introduce then a new  approach based on bid and ask prices observed for options. This last approach only takes  informations from the market, and leads to reduced spreads compatible with any observed spread that one wants to take into account.

\section{Bid-Ask Dynamic Pricing Procedure}
\label{sec:pricingfunction} 

\subsection{The economic model}
\label{sec:pricingfunction_model} 
  
Throughout this paper we work with  a filtered probability space denoted \ $(\Omega,{\cal F}_{\infty},({\cal F}_t)_{ t \in \R^+}, P)$. The 
filtration $({\cal F}_t)_{ t \in \R^+}$ satisfies the usual assumptions of right continuity and completedness and ${\cal F}_0$ is assumed to be the $\sigma$-algebra generated by the $P$ null sets of ${\cal F}_{\infty}$.  We assume that the time horizon is infinite, which is the most general case. Indeed if the time horizon is finite equal to $T$ we define ${\cal F}_s={\cal F}_T$ for any $s \geq T$. This general setting covers also the discrete time case. \\ 
 We want to construct a global dynamic pricing procedure which assigns to any financial product (asset, bond, option, basket, portfolio...) a dynamic ask (resp. bid) price process. 
 We assume that there is a  non risky liquid asset $S^0$ and  we  take it as num\'eraire. So from now on,
$(S^0)_t=1 \; \forall \; t \in \R^+$. Any financial asset is expressed in this num\'eraire.
\\ Due to the importance of stopping times in finance we will work with stopping times with values into $[0,+\infty]$.  
For any stopping time $\tau$, let ${\cal F}_{\tau}$ be the $\sigma$-algebra defined by
${\cal F}_{\tau}=\{A \in {\cal F}_{\infty}|\;\forall \; t \in \R^+\; A \INTER \{\tau \leq t \} \in {\cal F}_t\}$. 
Denote $L^{\infty}(\Omega,{\cal F}_{\tau},P)$ the Banach algebra of essentially bounded  real valued ${\cal F}_{\tau}$-measurable functions. We will always identify an essentially bounded \ ${\cal F}_{\tau}$-measurable function with its class in $L^{\infty}(\Omega,{\cal F}_{\tau},P)$.
In what follows a financial position at a stopping time $\tau$ means an element of $L^{\infty}({\cal F}_{\tau})$.
 Accordingly with previous literature dealing with surreplication prices and also with the paper of Jouini and Kallal we  modelize the dynamic ask price (or dynamic buying price) of $X$.
  As in \cite{JK}, \cite{BN02}, \cite{S} and \cite{Peng06}, we consider that selling $X$ is the same as buying $-X$. Therefore  the dynamic ask price of  $X$ is modeled as $\Pi_{\sigma,\tau}(X)$,(where $\sigma$ is a  stopping time $\sigma \leq \tau$)  and its dynamic bid price as $-\Pi_{\sigma,\tau}(-X)$. 
As  the process describing the non risky asset is constant, this leads to make the hypothesis of
 normalization $\Pi_{\sigma,\tau}(0)=0$, and of translation invariance for $\Pi_{\sigma,\tau}$ (we refer also to \cite{FG} and \cite{KS} for a discussion of this hypothesis).  The observation of the book of order leads to the conclusion that for $\lambda$ large enough, $\Pi_{0,\tau}(\lambda X)> \lambda \Pi_{0,\tau}(X)$ if $X$ is not perfectly liquid. Also taking into account the impact of diversification, it is natural to impose that  the ask price $\Pi_{\sigma,\tau}(X)$ is a  convex function of $X$.\\
Moreover for $\nu \leq \sigma \leq \tau$ it is natural to ask that we get the same result for the ask (resp bid) price at   time $\nu$ of a bounded ${\cal F}_{\tau}$-measurable financial instrument  if we evaluate it   either directly as $\Pi_{\nu,\tau}(X)$  or indirectly  as $\Pi_{\nu,\sigma}(\Pi_{\sigma,\tau}(X))$. This is the time consistency condition.\\
This motivates the following definition:

 \begin{definition}
 A dynamic pricing procedure  $(\Pi_{\sigma, \tau})_{0 \leq \sigma \leq \tau}$ on the filtered probability space $(\Omega,{\cal F}_{\infty},({\cal F}_t)_{t \in \R^+},P)$ (where $\sigma \leq \tau$ are  stopping times) is a 
family of maps $(\Pi_{\sigma, \tau})_{0 \leq \sigma \leq \tau}$, 
defined on $L^{\infty}({\cal F}_{\tau})$ with values into $L^{\infty}({\cal F}_{\sigma})$ satisfying the following 
four properties: 
\ben 
\item  monotonicity: $$\forall \; (X,Y) \in (L^{\infty}({\cal F}_{\tau}))^2, \;\;\;\mbox{if }\; X\leq Y \;\;\mbox{then } \;\;
\Pi_{\sigma,\tau}(X) \leq \Pi_{\sigma,\tau}(Y)$$  
\item  translation invariance: $$\forall \; Z \in L^{\infty}({\cal F}_{\sigma})\;,\;\;\forall \; X \in L^{\infty}({\cal F}_{\tau}) 
\;\;\Pi_{\sigma,\tau}(X+Z)=\Pi_{\sigma,\tau}(X)+Z$$
\item  convexity: $$\forall \; (X,Y) \in (L^{\infty}({\cal F}_{\tau}))^2 \;\;\forall \; \lambda \in [0,1]$$
$$\Pi_{\sigma,\tau}(\lambda X +(1-\lambda)Y) \leq \lambda \Pi_{\sigma,\tau}(X)+(1-\lambda) \Pi_{\sigma,\tau}(Y)$$
\item  normalization:  $\Pi_{\sigma,\tau}(0)=0$\\
\een
 For any $X \in L^{\infty}({\cal F}_{\tau})$, the dynamic ask (resp. bid) price process of $X$ is $\Pi_{\sigma,\tau}(X)$ (resp. $-\Pi_{\sigma,\tau}(-X)$).\\

\label{definition1}
\end{definition}
Notice first that as consequence of convexity and of normalization, the dynamic ask price $\Pi_{\sigma,\tau}(X)$ is always greater or equal to the dynamic bid price $-\Pi_{\sigma,\tau}(-X)$ (indeed $0 \leq \frac{1}{2}[ \Pi_{\sigma,\tau }(X) +\Pi_{\sigma,\tau}(-X)]$).\\
 The above  definition is very close to  the definition of non linear  expectations of Peng \cite{Peng} (with stopping times instead of deterministic times).
It is also closely related to the definition of conditional monetary utility functionals \cite{CDK} which was given (with stopping times) for processes in a discrete time setting. Another definition of pricing process $\Phi_t(X)$ with just one deterministic parameter can be found in Kl\"oppel and Schweizer \cite{KS}: it  corresponds to the bid price process $-\Pi_{t,\infty}(-X)$. 
 Remark also that $\rho_{\sigma,\tau}(X)=\Pi_{\sigma,\tau}(-X)$ is a  normalized dynamic risk measure as defined in \cite{BN04} (i.e. for any $\sigma \leq \tau$, $\rho_{\sigma,\tau}$ is a normalized risk measure on $L^{\infty}({\cal F}_{\tau})$ conditional to $L^{\infty}({\cal F}_{\sigma})$).\\
This definition of pricing procedure which does not assume homogeneity (for $\lambda>1,\;\;\Pi_{\sigma,\tau}(\lambda X) \geq \lambda \Pi_{\sigma,\tau}(X)$ (but not equal in general) takes into account the risk of liquidity. Also it  assigns to any financial product a dynamic bid price process and a dynamic ask price process. As in  \cite{JK} and \cite{J2} the bid ask spread reflects the transaction costs but also the real bid ask spread observed in the book of orders. Also as  noticed in the empirical study of Hamon and Jacquillat \cite{HJ}, the relative bid ask spread associated with one given asset is not constant during the same day. Therefore the relative bid ask spread  cannot be constant as it would be the case if it would be only the consequence of proportional transaction costs.
\begin{definition}
 A dynamic pricing procedure is called sublinear if 
$$\forall \lambda>0\;\; \forall X \in  L^{\infty}({\cal F}_{\tau})\;\; 
\Pi_{\sigma,\tau}(\lambda X)= \lambda \Pi_{\sigma,\tau}(X)$$
\end{definition}
This terminology is justified because convexity and homogeneity imply sublinearity for the ask price process.

\begin{definition}
 A dynamic pricing procedure  $(\Pi_{\sigma,\tau})_{0 \leq \sigma \leq \tau}$
is  time-consistent if 
$$\forall \; 0 \leq \nu \leq \sigma \leq \tau \;\; \forall \; X \in L^{\infty}({\cal F}_{\tau})\;\; 
\Pi_{\nu,\sigma}(\Pi_{\sigma,\tau}(X))=\Pi_{\nu,\tau}(X).$$ 

\label{def2}
\end{definition}
This definition of time consistency is a version with stopping times of the definition  with deterministic times first given by Peng  \cite{Peng}.
As mentioned above this property, which means that the price at time $\nu$ of any financial product defined at time $\tau \geq \nu$ can be indifferently computed either directely or using an intermediate instant of time,  is very natural in a context of pricing. The condition of time consistency expressed for stopping times is stronger than if it would just stated for deterministic times.
Notice also that for a time-consistent dynamic pricing procedure, the time consistency, the normalization and the translation invariance imply that  $\Pi_{\sigma,\tau}$ is the restriction of $\Pi_{\sigma,\infty}$ to $L^{\infty}({\cal F}_{\tau})$.  
Key  tools in the study of time-consistent dynamic pricing procedures  are the existence of a dual representation in terms of probability measures and the characterization of time consistency in terms of cocycle condition. 
The existence of  such a dual representation is  equivalent to continuity from below (cf. \cite{FS04} and \cite{DS}). 

\begin{definition} 
A dynamic pricing procedure $(\Pi_{\sigma, \tau})_{0 \leq \sigma \leq \tau}$ is continuous from below (resp.  above) 
if for any increasing (resp. decreasing)
sequence $X_n$ of elements of  $L^{\infty}({\cal F}_{\tau})$ such that $X=\lim\;X_n$, 
the increasing (resp. decreasing) 
sequence $\Pi_{\sigma,\tau}(X_n) $ has the limit $\Pi_{\sigma,\tau}(X)$.\\

\label{definition4}
\end{definition}

\begin{remark}
1. Continuity from below (resp above) for the dynamic pricing procedure means continuity from above (resp below) for the corresponding dynamic risk measure 
$\rho_{\sigma,\tau}(X)=\Pi_{\sigma,\tau}(-X)$\\
2. Continuity from above implies continuity from below (cf \cite{FS04} and \cite{DS}).\\
3. In \cite{JK}, lower semi continuity for the ask price process is assumed (i.e. if a sequence $X_n$ tends to $X$ and 
if $\Pi_{\sigma,\tau}(X_n)$ tends to $c$ then $\Pi_{\sigma,\tau}(X)\leq c$. It is then easy to verify that this implies the continuity from below for the dynamic pricing procedure. 
\end{remark}
 \begin{notation}
In all the following a TCPP means a time-consistent dynamic pricing procedure continuous from below.
\end{notation}

\begin{remark} A TCPP assigns a dynamic ask price process $\Pi_{\sigma, \infty}(X)$ to any european bounded option and also to any path dependent bounded option and also to any bounded
option defined on a basket, as any of these instruments can be represented by a bounded ${\cal F}_{\infty}$-measurable function.\\ When $X$ is no more bounded we can define as in \cite{CDK2}, $\Pi_{\sigma,\infty}(X)$ by: 
$$\Pi_{\sigma,\infty}(X)=lim_{n \rightarrow \infty}(lim_{m \rightarrow -\infty}(\Pi_{\sigma,\infty}(sup(inf(X,n),m)))$$ 
 For an American option of maturity time $\tau$, we can define $$\Pi_{\nu,\tau}(X)=\esssup_{\nu \leq \sigma \leq \tau}\Pi_{\nu,\sigma}(Y_{\sigma})$$ where $Y_{\sigma}$ is the corresponding European option of maturity $\sigma$.
\end{remark}
As already noticed, for any TCPP $(\Pi_{\sigma, \tau})_{0 \leq \sigma \leq \tau}$ $\rho_{\sigma,\tau}(X)=\Pi_{\sigma, \tau}(-X)$ is a time consistent dynamic risk measure continuous from above. Therefore we recall in that setting the two fundamental results which are the existence of a dual representation and the characterization of time consistency in terms of cocycle conditions which are both key properties for the study of TCPP.   
\subsection{General Dual representation }
In this subsection, we recall known  dual representation results.\\
Recall that from the study of Detlefsen and Scandolo \cite{DS} on conditional risk measures, on a probability space $(\Omega,{\cal F},P)$, it follows that any dynamic pricing procedure on \ $(\Omega,{\cal F}_{\infty},({\cal F}_t)_{ t \in \R^+}, P)$ continuous from below has a dual representation in terms of probability measures absolutely continuous with respect to $P$.
 Such a dual representation has been generalized in \cite{BN04} for any probability measure $Q \ll P$, to the projection on $L^{\infty}(\Omega,{\cal F}_{\sigma},Q)$. 
 For any pricing procedure continuous from below, for any probability measure $Q \ll P$, the projection of $\Pi_{\sigma,\tau}$ on $L^{\infty}(\Omega,{\cal F}_{\sigma},Q)$ has the following dual representation:
\begin{equation}
\forall \; X \in L^{\infty}({\cal F}_{\tau}) \;\;\Pi_{\sigma,\tau}(X)= Q-\esssup_{R\in  {\cal M}^1_{\sigma,\tau}(Q)}(E_R(X|{\cal F}_{\sigma})-\alpha^m_{\sigma,\tau}(R))\;\;Q\;a.s.
\label{eq_3}    
\end{equation}
where
\begin{equation}
{\cal M}^1_{\sigma,\tau}(Q)=\{R\;\mbox{on}\;(\Omega,{\cal F}_{\tau})\;,\; R\ll P,\;R_{|{\cal F}_{\sigma}}=Q\; \mbox{and } \; E_R(\alpha^m_{\sigma,\tau}(R)) < \;\infty\} 
\label{eq_2}
\end{equation}
and for any $R\ll P$,\begin{equation}
\alpha^m_{\sigma,\tau}(R)  =  R -\esssup_{X \in L^{\infty}(\Omega, {\cal F}_{\tau},P)} 
(E_R(X| {\cal F}_{\sigma})-\Pi_{\sigma,\tau}(X)) 
\label{eq_1}
\end{equation}
 In the particular case where $Q=P$ this result can be found in \cite{DS}. However the dual representation result in its most general form will be needed in section \ref{sec:admissible}.

Notice that  it follows from the normalization condition  that for any $Q$, the minimal penalty $\alpha^m_{\sigma,\tau}(Q)$ is always non negative.\\

When $\sigma=0$, 
${\cal M}^1_{0,\tau}(P)=\{R\;on\;(\Omega,{\cal F}_{\tau}),\; R\ll P,\;\alpha^m_{0,\tau}(R) < \;\infty\}$; This set will also be  denoted by  ${\cal M}_{0,\tau}$, and
\begin{equation}
\Pi_{0,\tau}(X)= \sup_{Q \in  {\cal M}_{0,\tau}}(E_Q(X)-\alpha^m_{0,\tau}(Q))
\label{eq_4}    
\end{equation}
This dual representation for monetary risk measures, i.e. in static case, was proved by F\"ollmer and  Schied \cite{FS04}.

\begin{remark}
In case of continuity from above, the ``$\esssup$'' (resp. ``$\sup$'') in the above equations (\ref{eq_3}) (resp. (\ref{eq_4})) are in fact realized. (cf. \cite{FS02} for the static case and \cite{BN01} for the conditional case). Therefore in case of continuity from above, for any $\tau$, as $\Pi_{0,\tau}(0)=0$, there is always a probability measure $Q\ll P$ with $0$ minimal penalty, $\alpha^m_{0,\tau}(Q)=0$.
\label{rq1}
\end{remark}
Theorems of representation have been proved in the static case for monetary convex risk measures in \cite{FS02}, \cite{FS04} and \cite{FG}. In the conditional or dynamic case they have been proved in \cite{DS}, \cite{BN01}, \cite{KS} and \cite{BN04}. 

\subsection{Characterization of time consistency}
 Recall from \cite{BN04} the following  results characterizing time consistency:\\  A dynamic pricing procedure  $(\Pi_{\sigma,\tau})_{0 \leq \sigma \leq \tau}$ continuous from below,
 is time-consistent if and only if for any probability measure $Q$ absolutely continuous with respect to $P$, the minimal prenalty function satisfies the following cocycle condition for any stopping times $\nu \leq \sigma \leq \tau$:
\begin{equation}
\alpha^m_{\nu, \tau}(Q)=\alpha^m_{\nu, \sigma}(Q)+E_Q(\alpha^m_{\sigma, \tau}(Q)|{\cal F}_{\nu})\; \;Q\;a.s
\label{eq_5}
\end{equation}
The cocycle condition appeared for the first time in \cite{BN03}. The characterization of time consistency in terms of cocycle condition for the minimal penalty has been proved  in \cite{BN03} and  then in \cite{FP} with slightly different hypothesis, in both cases under restrictive conditions. The characterization in the general  case of dynamic risk measures continuous from above is proved in \cite{BN04}.\\
Another characterization of time consistency was given in Cheridito, Delbaen, Kupper \cite{CDK} in terms of acceptance sets and also in terms of a concatenation condition. The advantage of the cocycle condition is that it is easy to check.

\section{No Free Lunch Pricing Procedure}
\label {sec:No Free Lunch}

In this section we define two intrinsic notions of No Free Lunch for a TCPP: one static and one multiperiod. We then prove that due to the time consistency condition, these notions are in fact equivalent.
The main result of this section is  the Fundamental Theorem of pricing proved in that context. It generalizes   the result of Jouini and Kallal  \cite{JK}:
We prove that a TCPP has No Free Lunch if and only if there is an equivalent probability  measure that transforms some process between the bid and ask price process into a  martingale.
We can notice that the case studied in \cite{JK} corresponds in our setting to the case of sublinear dynamic pricing (i.e. to the case where there is no liquidity risk). However an important difference is that a TCPP is always defined on the set of all bounded  measurable functions. Notice that here  we do not assume a priori that the bid price process (and the ask price process) are  right continuous.\\The notion of No Free Lunch in our setting has an easy definition because any TCPP is defined on the set of all bounded measurable functions.
A static arbitrage is a financial instrument, i.e. a bounded ${\cal F}_{\infty}$-measurable  function $X$, non negative and non zero, that one can buy at time zero at a non positive price, that is $\Pi_{0,\infty}(X) \leq 0$. Equivalently there is No Static Arbitrage if and only if  ${\cal C}  \INTER L^{\infty}_+(\Omega,{\cal F}_{\infty},P)=\{0\}$
where   
$${\cal C}=\{\lambda X, (\lambda,X) \in  \R^+ \times L^{\infty}(\Omega,{\cal F}_{\infty},P)\;/\; \Pi_{0,\infty}(X) \leq 0\}$$
 As usual this definition of no arbitrage is not sufficient, we have to pass to the notion of No Free Lunch. 
\begin{definition}
The Dynamic Pricing Procedure admits No Static Free Lunch if 
$\overline {\cal C}  \INTER L^{\infty}_+(\Omega,{\cal F}_{\infty},P)=\{0\}$
where $\overline {\cal C}$ is the weak* closure of  ${\cal C}$ (i.e. the closure of ${\cal C}$ in $ L^{\infty}(\Omega,{\cal F}_{\infty},P)$ for the topology $\sigma(L^{\infty}(P),L^1(P))$) 
$${\cal C}=\{\lambda X, (\lambda,X) \in  \R^+ \times  L^{\infty}(\Omega,{\cal F}_{\infty},P)\;/\; \Pi_{0,\infty}(X) \leq 0\}$$

\label{defStatic}
\end{definition}
   
 We adapt now to our framework the definition of no multiperiod free lunch given by Jouini and Kallal \cite{JK} in the case of sublinear pricing. In particular here the time horizon is infinite, the pricing procedure associates to any bounded ${\cal F}_{\infty}$-measurable function a dynamic bid (resp ask) price process, and the dynamic pricing procedure is not homogeneous. This motivates the following definition:
\begin{definition}
The set of financial products atteignable at zero cost via simple strategies is \begin{eqnarray}
{\cal K}_0=\{X=X_0+\sum_{1 \leq i \leq n}(Z_i-Y_i)\;/\;\Pi_{0,\infty}(X_0) \leq 0;\nonumber\\
 \Pi_{\tau_i,\infty}(Z_i) \leq -\Pi_{\tau_i,\infty}(-Y_i) \;\forall \;1 \leq i \leq n\}\nonumber
\end{eqnarray}
where $0 \leq \tau_1 \leq ...\leq \tau_n < \infty$ are stopping times.
\label{Zerocost}
\end{definition}

This means that one can get $X_0$ at time zero at zero or negative cost and then at time $\tau_i$, one sells $Y_i$ and buys $Z_i$ in a self-financing way.
Of course in the preceding definition, due to the translation invariance property of $\Pi_{\tau_i,\infty}$, one can restrict to non negative $Y_i$
and $Z_i$. It doesn't change the set ${\cal K}_0$. As in the static case  we define as follows the No Mutiperiod Free Lunch condition.

\begin{definition} The dynamic pricing procedure  has No Multiperiod Free Lunch if  $\overline {\cal K}  \INTER L^{\infty}_+(\Omega,{\cal F}_{\infty},P)=\{0\}$
where $\overline {\cal K}$ is the weak* closure of 
$${\cal K}=\{\lambda X, (\lambda, X) \in \R^+ \times {\cal K}_0\}$$
\label {Multi}
\end{definition}
We can now prove our  main theorem:

\begin{theorem}
Let $(\Pi_{\sigma,\tau})_{0 \leq \sigma \leq \tau}$ be a TCPP. The following conditions are equivalent:\\
i) The TCPP has No Muliperiod Free Lunch\\
ii) The TCPP has No Static Free Lunch.\\
iii) There is a probability measure $R$ equivalent to $P$ with zero  minimal penalty: $\alpha^m_{0,\infty}(R)=0$\\
iv) There is a probability measure $R$ equivalent to $P$ such that for any stopping time $\sigma$,
\begin{equation}
\forall X \in  L^{\infty}(\Omega,{\cal F}_{\infty},P)\;\; \; -\Pi_{\sigma,\infty}(-X)\leq E_R(X|{\cal F}_{\sigma}) \leq \Pi_{\sigma,\infty}(X)
\label{eqwr}
\end{equation}

\label{thmarbitrage}
\end{theorem} 
{\bf Proof}:
- Start with the proof of {\it iv) implies {\it i)}}: 
  Let $X \in {\cal K}_0$. From equation (\ref{eqwr}) applied to $Y_i$ and $Z_i$, and the self financing condition, we get 
$$E_R(Z_i-Y_i|{\cal F}_{\tau_i}) \leq  \Pi_{\tau_i,\infty}(Z_i) + \Pi_{\tau_i,\infty}(-Y_i) \leq 0$$
Also from equation (\ref{eqwr}), $E_R(X_0) \leq 0$.
Thus $\forall X \in {\cal K}_0 \;\; E_R(X) \leq 0$.
As $R$ is a probability measure equivalent to $P$, $E_R$ is linear and continuous for the weak* topology  $\sigma(L^{\infty}(P),L^1(P))$ and
therefore  $\forall X \in \overline {\cal K}\;\; E_R(X) \leq 0$ and thus $\forall X \in \overline {\cal K}\INTER L^{\infty}_+(\Omega,{\cal F}_{\infty},P)$, $X=0 \;\; R\;a.s.$. $R$ is equivalent to $P$ so $X=0$, establishing the No Multiperiod Free Lunch.\\
- {\it i)} implies {\it ii)} is trivial (taking $Y_i=Z_i=0$ for any $i$).\\
- {\it ii)} implies {\it iii)}: 
The proof is in two steps as for the proof  of Kreps Yan theorem (cf Theorem 5.2.2. of \cite{DS02}).\\
- First step, Hahn Banach separation: 
As $\Pi_{0,\infty}$ is convex, $\overline {\cal C}$ is a  convex cone in $ L^{\infty}$ closed for the topology $\sigma( L^{\infty},L^1)$. For any $Y \in   L^{\infty}_+(\Omega,{\cal F}_{\infty},P )$, $Y \neq 0$, $\{Y\}$ is compact and  $\overline {\cal C}\INTER \{Y\}= \emptyset$. From Hahn Banach theorem, there  is a $g_Y \in L^1$ such that $E(g_Y Y)> \sup_{X \in \overline {\cal C}} E(g_YX)$. As $\overline {\cal C}$ is a cone, containing   $L^{\infty}_-$ (from monotonicity and normalization), it follows that $\sup_{X \in \overline {\cal C}} E(g_YX)=0$ and  $g_Y \geq 0$.\\
- Second step, exhaustion argument exactly as in the proof of Theorem 5.2.2. of \cite{DS02}. Consider ${\cal G}=\{g \in L^1+ , \;\sup_{X \in \overline {\cal C}} E(gX)=0\}$. Let $g_0 \in {\cal G}$ such that 
$$P(g_0>0)=\sup \{P(g>0), \;g \in {\cal G}\}$$
As in  \cite{DS02}, the probability measure $R$ of Radon Nikodym derivative $\frac{g_0}{E(g_0)}$ is equivalent to $P$ and $\sup_{X \in \overline {\cal C}} E_R(X)=0$.\\
$\forall X \in L^{\infty}(\Omega,{\cal F}_{\infty},P )$, from translation
invariance, $X-\Pi_{0,\infty}(X)$ belongs to ${\cal C}$ and therefore $E_R(X) \leq 
\Pi_{0,\infty}(X)$. Applying this inequality to $-X$, we get
\begin{equation}\forall X \in  L^{\infty}(\Omega,{\cal F}_{\infty},P )\;\;-\Pi_{0,\infty}(-X)
\leq E_R(X) \leq \Pi_{0,\infty}(X)
\label{eq_NFL}
\end{equation}
From the definition of the minimal penalty 
$$\alpha^m_{0,\infty}(R)=\sup_{X \in  L^{\infty}(\Omega,{\cal F}_{\infty},P )}(E_R(X)-\Pi_{0,\infty}(X))\leq 0$$ 
On the other hand as already noticed, the minimal penalty is non negative so $\alpha^m_{0,\infty}(R)=0$ and {\it iii)}
is satified.\\
-{\it iii)} implies {\it iv)}: From the cocycle condition and the non negativity of the minimal penalty, it follows that $\alpha^m_{\sigma,\infty}(R)=0\;\; \forall \sigma $.  Then {\it iv)} follows from the dual representation of $\Pi_{\sigma,\infty}$, equation (\ref{eq_3}), as $R$ is equivalent to $P$. (Notice that for this proof we need only the usual  dual representation and not its extended shape). 
\begin{corollary} Any No Free Lunch TCPP has a representation in terms of equivalent probability measures, i.e.
\begin{equation}
\forall \; X \in L^{\infty}({\cal F}_{\tau}) \;\;\Pi_{\sigma,\tau}(X)= \esssup_{Q \in  {\cal M}^{1,e}_{\sigma,\tau}(P)}(E_Q(X|{\cal F}_{\sigma})-\alpha^m_{\sigma,\tau}(Q))
\label{eq_NFL2}    
\end{equation}
where
\begin{equation}
{\cal M}^{1,e}_{\sigma,\tau}(P)=\{Q\;\mbox{on}\;(\Omega,{\cal F}_{\tau})\;,\; Q\sim  P,\;Q_{|{\cal F}_{\sigma}}=P\; \mbox{and } \; E_Q(\alpha^m_{\sigma,\tau}(Q)) < \;\infty\} 
\label{eq_2NFL}
\end{equation}
\label{corNFL}
\end{corollary}
{\bf Proof}:  
  From the preceding theorem  if the TCPP has No Free Lunch there is a probability measure equivalent to $P$ with zero penalty and then from Theorem 4 of \cite{KS} it has a representation in terms of equivalent probability measures.\\

\begin{proposition}

Let $(\Pi_{\sigma,\tau})_{\sigma \leq \tau}$ be a No Free Lunch  TCPP.
  For any $X \in L^{\infty}(\Omega,{\cal F}_{\infty},P)$,  for any probability measure $R$ equivalent to $P$ with zero penalty, $(\Pi_{\sigma,\infty}(X))_{\sigma}$  is a  $R$-supermartingale process with a cadlag modification.
(resp. $-(\Pi_{\sigma,\infty}(-X))_{\sigma}$ is  a $R$-submartingale process with a cadlag modification), and 
\begin{equation}
-\Pi_{\sigma,\infty}(-X)) \leq E_{R}(X|{\cal F}_{\sigma}) \leq \Pi_{\sigma,\infty}(X) 
\label{eq01}
\end{equation}
 for any stopping time $\sigma$.
\label{propNFL}
\end{proposition}
{\bf Proof}:  From the preceding theorem, there is a probability measure $R$ equivalent to $P$ with zero penalty. The result follows then  from Theorem 2.3. of \cite{BN04}.

 Notice that the regularity of the paths of the ask (resp. bid) price processes is obtained here as a consequence of the No Free Lunch condition and the time consistency. On the contrary in  \cite{JK} the regularity of paths was assumed.\\
 Using the preceding theorem we can characterize the No Free Lunch sublinear TCPP.
\begin{corollary}
For any No Free Lunch sublinear TCPP  there is a stable subset ${\cal Q}$ of the set of probability measures equivalent to $P$  such that for any stopping times $\sigma \leq \tau$, for any $X \in L^{\infty}({\cal F}_{\tau})$,
\begin{equation}
\Pi_{\sigma,\tau}(X)=\esssup_{Q \in {\cal Q}}E_Q(X|{\cal F}_{\sigma})
\label{eqsub0}
\end{equation}
Conversely the equation (\ref{eqsub0}) defines a No Free Lunch sublinear TCPP. 
\label{corNFLS}
\end{corollary}
{\bf Proof}: From Corollary \ref{corNFL}, a No free Lunch TCPP has a representation in terms of equivalent probability measures. As a sublinear pricing procedure is homogeneous, any probability measure $Q$ with finite penalty $\alpha^m_{0,\infty}(Q)<\infty$ has in fact a zero penalty. Denote ${\cal Q}=\{Q \sim  P \;| \;\alpha^m_{0,\infty}(Q)=0\}$. For any stopping times $\sigma \leq \tau$, and any $X \in L^{\infty}({\cal F}_{\tau})$, $$\Pi_{\sigma,\tau}(X)=\esssup_{Q \in {\cal Q}}E_Q(X|{\cal F}_{\sigma})$$  
From Delbaen \cite{D}, the time consistency is then equivalent to the stability of  ${\cal Q}$.\\
Conversely,  the dynamic sublinear pricing procedure defined by the formula  (\ref{eqsub0}) is continuous from above. It is time consistent from \cite{D}. It has No Free Lunch from Theorem \ref{thmarbitrage}.

The following corollary gives a  sufficient condition for a pricing procedure
 to have No Free Lunch. Recall the definition of non degeneracy \cite{BN04}.

\begin{definition}
A dynamic pricing procedure is non degenerate if   $$\forall \; A \in {\cal F}_{\infty}\;\;  [\Pi_{0,\infty}(\lambda 1_A)=0\;\;\forall \; \lambda \in {\R_{+}}^*]\;\; implies \;\;P(A)=0$$ 
\label{defnd}
\end{definition}

\begin{corollary}
Let $(\Pi_{\sigma,\tau})_{\sigma \leq \tau}$ be a  non degenerate TCPP. Assume that  there is a probability measure $Q$ absolutely continuous with respect to $P$ with zero minimal penalty (i.e. $\alpha^m_{0,\infty}(Q)=0)$.  Then this TCPP has No Free Lunch. In particular:\\ 
Any non degenerate  sublinear TCPP has No Free Lunch.\\
 Any non degenerate TCPP  continuous from above has No Free Lunch.
\label{cor ND}
\end{corollary}
{\bf Proof}: Assume that  $(\Pi_{\sigma,\tau})_{\sigma \leq \tau}$ is non degenerate. From  Lemma 2.4. of \cite{BN04} any probability measure $Q \ll P$ with zero minimal penalty is equivalent to $P$. When the pricing procedure is sublinear any probability measure with finite minimal penalty has zero minimal penalty. Thus there is at least one  $Q$ equivalent with  $P$ with zero minimal penalty. And from Remark \ref{rq1}, in case of continuity from above  there is always  a probability measure  $Q \ll P$ with zero minimal penalty. $Q$ is equivalent with  $P$ as $(\Pi_{\sigma,\tau})_{\sigma \leq \tau}$ is non degenerate.
The result is then a consequence of Theorem \ref{thmarbitrage}.

One  ends this section with a corollary which will be very usefull  for the construction of No Free Lunch TCPP. One can notice that in this corollary the conditions required for the penalty are only the cocycle condition and the locality. In many examples it is easy to verify if these conditions are satisfied. On the contrary one can insist on the fact  that one does not need to know if the penalty is the minimal one.
\begin{corollary}
Let ${\cal Q}$ be a  stable set of probability measures on $(\Omega,{\cal F}_{\infty})$ all equivalent to $P$. let $\alpha$ be a non negative local penalty function defined on ${\cal Q}$, satisfying the cocycle condition . Assume that there is at least one probability measure $Q \in {\cal Q}$ with zero penalty (i.e. $\alpha_{0,\infty}(Q)=0$). Then $(\Pi_{\sigma,\tau})_{0 \leq \sigma \leq \tau}$  defined by
\begin{equation}
\Pi_{\sigma,\tau}(X)=\ esssup _{Q \in {\cal Q}} \{E_Q(X|{\cal F}_{\sigma}))-\alpha_{\sigma,\tau}(Q)\} 
\label{eqdef}
\end{equation}
is a No Free Lunch TCPP.
\label{proprecall}
\end{corollary}
{\bf Proof}: Equation (\ref{eqdef}) defines a pricing procedure continuous from above. From \cite{BN03} in case of deterministic times and \cite{BN04}
for the extension to stopping times, it follows that this pricing procedure is time consistent. From the minimality of $\alpha^m_{0,\infty}$, any probability measure $Q$ such that $\alpha_{0,\infty}(Q)=0$ satisfies also $\alpha^m_{0,\infty}(Q)=0$. Thus from Theorem \ref{thmarbitrage}, the TCPP has No Free Lunch.

Up to now, we have studied TCPP in a model free setting. We have characterized the intrinsic notion of
No Free Lunch for a TCPP by the existence of an equivalent probability measure
that  transforms a process between the  bid process and the ask process of any financial instrument into a martingale.
In the next Section, we make use of TCPP to reduce the bid ask spread in the more usual setting
where the dynamics of some reference assets are assumed to be known.

\section{Reduction of the bid ask spread}
\label {sec:admissible}  

\subsection{The market model}
In this section we consider the following general and classical framework of financial markets: Let  
$(\Omega,{\cal F}_{\infty},({\cal F}_t)_{ t \in \R^+}, P)$
be a filtered probability space satisfying the usual assumptions.  
Let $(S^k)_{0 \leq k \leq d}$, be a 
$\R^{d+1}$ valued adapted stochastic process. The problem of pricing in this framework is the problem of providing a dynamic pricing procedure 
which extends the dynamics of the reference assets $(S^k)_{0 \leq k \leq d}$. 
Any linear dynamic pricing procedure is given by the conditional expectation with respect to an equivalent 
local martingale measure for $(S^k)_{0 \leq k \leq d}$. Therefore the classical ask price is defined as the surreplication price $\sup_{Q \in {\cal M}^e}E_Q(X)$ where ${\cal M}^e$ denotes the set of all equivalent local martingale measures. It is well known that the corresponding bid-ask spread is too wide and therefore different methods have been developed in the literature in order to obtain sharper bounds. One of them is the so called ``No good deal'' pricing. Another method consists in restricting the set of admissible strategies. Both of them can be written in our general framework.\\
  As before, we assume that there is a liquid positive asset $S^0$ and we take it as num\'eraire. The other reference assets are modelized by locally bounded adapted  processes $(S^k)_{1 \leq k \leq d}$, i.e. there exists a sequence $\tau_n$ of stopping times increasing to $\infty$ such that the stopped processes $(S^k)^{\tau_n}_t=S^k_{t \wedge \tau_n}$ are uniformly bounded for each $n$. This assumption is not too restrictive, indeed it is satisfied for any cadlag (right continuous with left limit) process with uniformly bounded jumps.
Notice that in mathematical finance it is important to allow for processes with jumps and also not to restrict to the filtration generated by a multidimensional Brownian motion.\\We start with some general considerations on TCPP extending the dynamics of the reference assets $(S^k)$.

\begin{definition}The TCPP extends the dynamics of the process $(S^k)_{0 \leq k \leq d}$ 
if $$\forall \; 0 \leq k \leq d,\;\;\forall \; n \in \Z \;  \;\forall \; 0 \leq \sigma \leq \tau \;\;such\;\; that \;\;S^k_{\tau} \in L^{\infty}({\cal F}_{\tau}), $$ 
$$\Pi_{\sigma,\tau}(n S^k_{\tau})= n S^k_{\sigma}\;\;(A1)$$
\end{definition}
  
In this definition we consider only integer multiples of the processes $S^k$ because there are the only one that can be traded.\\
We prove now that a TCPP continuous from below extends the dynamics of the process $(S^k)_{0 \leq k \leq d}$ if and only if it has a dual representation in terms of  local martingale measures for the process $(S^k)_{0 \leq k \leq d}$:
\begin{theorem}
A TCPP extends the dynamics of the process  $(S^k)_{0 \leq k \leq d}$ if and only if any probability measure $R \ll P$ such that $\alpha^m_{0,\infty}(R) < \infty$ is a local martingale measure with respect to the process  $(S^k)_{0 \leq k \leq d}$
i.e.
\begin{eqnarray}
\forall 0 \leq k \leq d\;\;\forall \tau \;\; such \;\;that \;\;S^k_{\tau}
\in L^{\infty}({\cal F}_{\tau}), \nonumber\\
\forall \sigma \leq \tau\;\;E_R(S^k_{\tau }|{\cal F}_{\sigma})=S^k_{\sigma}\;\;R\;a.s.
\label{eqLMM}
\end{eqnarray}
\label{propEMM}
\end{theorem}
{\bf Proof}:
1. Assume first that $(\Pi_{\sigma,\tau})_{0 \leq \sigma \leq \tau}$ extends the dynamics of the process $(S^k)_{0 \leq k \leq d}$.\\
 Let $R \ll P$ such that  $\alpha^m_{0,\infty}(R)< \infty$. From the cocycle condition, for any stopping times $\sigma \leq \tau$, $E_R(\alpha^m_{\sigma,\tau}(R))<\infty$.
For any positive integer $j$, denote $A_j=\{\omega\;/ \; \alpha^m_{\sigma,\tau}(R)(\omega) \leq j \;\}$ $A_j \in {\cal F}_{\sigma}$ and  $\;R((UA_j)^c)=0$. \\ Let $\tau$ and $k$ such that $S^k_{\tau} \in  L^{\infty}({\cal F}_{\tau})$. 
Then by hypothesis, $ \forall \; n \in \Z$, $\Pi_{\sigma,\tau}(n S^k_{\tau})= n S^k_{\sigma}$\\
For any $n \in \Z$, from the general dual representation equation (\ref{eq_3}),  
$$n S^k_{\sigma}=\Pi_{\sigma,\tau}(n S^k_{\tau}) \geq n E_R(S^k_{\tau} |{\cal F}_{\sigma})-\alpha^m_{\sigma,\tau}(R)\;R \;a.s.$$ 
As  $1_{A_j}\alpha^m_{\sigma,\tau}(R)$ is $R$ essentially bounded the restriction of the above inequality to $A_j$ for any $n \in\Z$, implies that 
 $E_R(S^k_{\tau} |{\cal F}_{\sigma})$ and $S^k_{\sigma}$ coincide $R$ a.s. on every $A_j$ and thus on $\Omega$.
This proves the equation (\ref{eqLMM}) for any  $ R $ such that $\alpha^m_{0,\infty}(R)< \infty$\\
 2.  Conversely assume that the local martingale property is satisfied. Let $k \in\{1,..d\}$.  Let ${\tau}$ be a stopping time such that $S^k_{\tau} \in L^{\infty}({\cal F}_{\tau})$. 
 Let $R \in {\cal M}^1_{\sigma,\tau}(P)$.\\ First step: We  prove that  
$E_{R}(S^k_{\tau}|{\cal F}_{\sigma})=S^k_{\sigma}$. Choose $Q_1 \in  {\cal M}_{0,\sigma}(P)$ and $Q_2 \in  {\cal M}^1_{\tau,\infty}(P)$. From the cocycle condition it follows that the probability measure $\tilde R$ with Radon Nikodym derivative $\frac{d \tilde R}{dP}=\frac{dQ_1}{dP} \frac{dR}{dP} \frac{dQ_2}{dP}$
satisfies $\alpha^m_{0,\infty}(\tilde R) < \infty$.  By definition of $\tilde R$,  
$$E_{\tilde R}(S^k_{\tau}|{\cal F}_{\sigma})=E_{R}(S^k_{\tau}|{\cal F}_{\sigma})$$
  By hypothesis, $\tilde R$ is a local martingale measure with respect to $S^k$
,
so $E_{R}(S^k_{\tau}|{\cal F}_{\sigma})=S^k_{\sigma}$. This ends the first step.\\
From the  dual representation (equation (\ref{eq_3})), applied with $Q=P$, we get
$$\Pi_{\sigma,\tau}(n S^k_{\tau})= (nS^k_{\sigma}+ {\esssup}_{R \in {\cal M}^1_{\sigma,\tau}(P)}(-\alpha^m_{\sigma,\tau}(R))=n S^k_{\sigma}+\Pi_{\sigma,\tau}(0)=n S^k_{\sigma}$$
 Thus the TCPP extends the dynamics of the process $(S^k)_{0 \leq k \leq d}$\\
From  Theorem \ref{propEMM} and Corollary \ref{corNFL} we deduce the following corollary:

\begin{corollary} A No Free Lunch TCPP $(\Pi_{\sigma \tau})_{0 \leq \sigma \leq \tau}$ extends the dynamics of the process $(S^k)_{0 \leq k \leq d}$ if and only if it has a dual representation in terms of equivalent local martingale measures for the  process $(S^k)_{0 \leq k \leq d}$. 
\label{corEMM}
\end{corollary}

From the non negativity of the minimal penalty, we then deduce the following result:
\begin{corollary}
The ask price associated with a No Free Lunch TCPP\\ $(\Pi_{\sigma \tau})_{0 \leq \sigma \leq \tau}$ extending the dynamics of the process $(S^k)_{0 \leq k \leq d}$
is never greater than the surreplication price. More precisely for any  $X$ in $L^{\infty}({\cal F}_{\tau})$,
\begin{equation}
\inf_{Q \in {\cal M}^e}E_Q(X) \leq  -\Pi_{\sigma \tau}(-X) \leq \Pi_{\sigma \tau}( X) \leq \sup_{Q \in {\cal M}^e}E_Q(X)
\label{eqbounds}
\end{equation}

Where ${\cal M}^e$ denotes the set of all equivalent local martingale measures with respect to the process $(S^k)_{0 \leq k \leq d}$.
\label{sur}
\end{corollary}
The bounds obtained for any No Free Lunch TCPP extending the dynamics of the process $(S^k)$ extend to general (convex) pricing procedures the usual bounds known for linear or sublinear pricing procedures. We can notice that these bounds have been obtained from the time consistency condition without using replicating or surreplicating portfolio.\\
We also get the following dual representation result for sublinear TCPP:

\begin{corollary}
For any No Free Lunch sublinear TCPP extending the dynamics of $(S^k)_{0 \leq k \leq d}$ there is a stable subset ${\cal Q}$ of the set ${\cal M}^e$ of equivalent local martingale measures for  $(S^k)_{0 \leq k \leq d}$ such that for any stopping times $\sigma \leq \tau$, for any $X \in L^{\infty}({\cal F}_{\tau})$,
\begin{equation}
\Pi_{\sigma,\tau}(X)=\esssup_{Q \in {\cal Q}}E_Q(X|{\cal F}_{\sigma})
\label{eqsub}
\end{equation}
Conversely the equation (\ref{eqsub}) defines a No Free Lunch TCPP extending the dynamics of $(S^k)_{0 \leq k \leq d}$.
\label{sublinearEMM}
\end{corollary}
{\bf Proof} It follows from Corollary \ref{corNFLS} and Theorem \ref{propEMM}.

Now, we consider different approaches in order to reduce the bid ask spread.
\subsection{Sublinear pricing and good deal bonds}
``No good deal'' pricing started with two papers published in 2000. One by Cochrane and Saa Requejo \cite{CSR}, the other by Bernardo and Ledoit \cite{BL}.
In both cases the idea is to exclude not only financial products leading to arbitrage but also products which are too attractive. The condition is then expressed in terms of a condition on the Radon Nikodym derivative $\frac{dQ}{dP}$ of the probability measures $Q$ that one can use to price. The ask price is then defined as $\sup_{Q \in {\cal Q}}E_Q(X)$ ( and the bid price as $\inf_{Q \in {\cal Q}}E_Q(X)$.
The static case is considered for example in \cite{CSR}, \cite{BL} and \cite{S}. Jaschke and K\"uchler have established the link between coherent static risk measures and no good deal pricing.\\
In a multiperiod setting, Cochrane and Saa Requejo \cite{CSR} consider a set of probability measures $Q$ defined  in the following way:
between instants of time $t_i$ and $t_{i+1}$, denote $m_i=\frac{E(\frac{dQ}{dP}|{\cal F}_{t_{i+1}})}{E(\frac{dQ}{dP}|{\cal F}_{t_{i}})}$. The condition on $Q$ is then of two kinds: the first one is that it is a martingale measure with respect to some basic assets and the second one is  $E(m_i^2|{\cal F}_{t_i})\leq A_i^2$. where $A_i$ are some given  ${\cal F}_{t_i}$-measurable functions. These two conditions clearly define  a stable set of probability measures ${\cal Q}$, leading thus to a sublinear TCPP.
In a dynamic setting, Cochrane and Saa Requejo consider the case where the processes describing the basic assets $S^k$ are driven by a multidimensional Wiener process. No good deal pricing has been extended to models with jumps by Bj\"ork and Slinko \cite{BS} and also by Kl\"oppel and Schweizer\cite{KS}, providing in both cases a sublinear pricing procedure. The ask price is given by $\sup_{Q \in {\cal Q}}E_Q(X)$. From Corollary \ref{corNFLS}, the pricing procedure is time consistent if and only if the set ${\cal Q}$ is stable.  Kl\"oppel and Schweizer have verified the stability of the set of equivalent probability measures associated with  their pricing procedure. This gives thus an example of a No Free Lunch  sublinear TCPP in the setting of Levy processes.

\subsection{TCPP from portfolio constraints}
\label {sec:admissible portfolio}
 Several authors  consider also the problem of pricing and hedging under constraints on the admissible porfolios. \\
- In  \cite{CK1}  Cvitanic and Karatzas  consider a financial market where portfolios are constrained to take values in a given closed convex set. The framework is that of diffusions driven by a multi dimensional Brownian motion.\\ 
-  F\"ollmer and Kramkov have extended the results of \cite{CK1} to the following general framework:  The discounted processes of the  basic assets are described by a locally bounded semimartingale $(S^k)_{0 \leq k \leq d}$. Under conditions on the set ${\cal H}$ of admissible strategies, they have computed the shape of the value process associated with the minimal ${\cal H}$-constrained  hedging portfolio. Under a slight modification,  Kl\"oppel and Schweizer have proved \cite{KS}, Proposition 10 that the value process associated with any $X \in L^{\infty}({\cal F}_T)$ is 
\begin{equation}
\Pi_{\tau,T}(X)=\esssup_{Q \in {\cal P}({\cal H})}(E_Q(X|{\cal F}_{\tau})-E_Q[{\cal A}^{\cal H}(Q)_T-{\cal A}^{\cal H}(Q)_{\tau}|{\cal F}_{\tau}])
\label{eqPC}
\end{equation} 
where ${\cal A}^{\cal H}(Q)$ is the smallest increasing predictable process $A$
such that  $Y-A$ is a local $Q$ supermartingale  for any $Y \in \{H.S\;|\; H \in {\cal H}\}$. From this definition of  ${\cal A}^{\cal H}(Q)$ it is not difficult to verify that the penalty $\alpha_{\tau,T}(Q)=E_Q[{\cal A}^{\cal H}(Q)_T-{\cal A}^{\cal H}(Q)_{\tau}|{\cal F}_{\tau}]$ is local,  that the cocycle condition (for stopping times) is satisfied and that ${\cal P}({\cal H})=\{Q \;|\;\alpha_{0,\infty}(Q)<\infty\}$ is stable. Therefore from Corollary \ref{proprecall} the formula (\ref{eqPC}) defines a No Free Lunch TCPP. Notice that these sufficiant conditions for time consistency are easy to verify from the formula (\ref{eqPC}).  The time consistency condition for deterministic times is  also verified by Kl\"oppel and Schweizer \cite{KS2}, in a different way,  using minimal strategies.

\subsection{Calibration of TCPP with observed bid ask prices for options}
\label{sec:calibration}
 We propose now a new approach in order to construct  bid ask spreads consistent with the bid ask spreads observed in a real financial market for a family of options. This approach   takes fully advantage of our axiomatic approach. In  financial markets the market-maker observes the bid and ask prices associated with some financial assets as options. Some of these options (vanilla options or even some exotic options) are quite liquid and therefore it is relevant to take them into account in order to provide a global TCPP. As pointed out by Avellaneda, Levy and Paras \cite{ALP}, the prices of options provide important information on the volatility. Such a pricing procedure, compatible with observed bid and ask prices for the reference options, can then be used in order to provide a bid and an ask price associated with some less liquid financial product. 
Therefore the aim of this section is to define TCPP not only   extending  the dynamics of the basic assets but also compatible with  observed bid and ask prices for  options. Notice that if one knows the process describing the dynamics of  basic assets (i.e. one knows the process $(S^k)_{0 \leq k \leq d}$), an option on one or several of these underlying assets is defined at its maturity date $T$ as a  ${\cal F}_T$-mesurable  function. However, the dynamic process associated with this option is not known.\\ The new conditions (compatibility with observed bid and ask prices for options) that we impose on the TCPP are all given by conditions obtained from the market. They are preference free. And we prove that they lead to  an upper bound for the ask price (resp. a lower bound for the bid price) which are better than the surreplication (resp subreplication ) price.   \\
 The market model is the following: Let $(\Omega,{\cal F}_{\infty},({\cal F}_t)_{ t \in \R^+}, P)$ be a filtered probability space. Consider a  reference family  $((S^k)_{0 \leq k \leq d},(Y^l)_{1 \leq l \leq p})$  composed of two kinds of assets:\\
- basic assets as in the usual setting, and among them an asset $S^0$ that we take as num\'eraire. The other   basic  assets are modeled by a locally bounded adapted  process $(S^k)_{1 \leq k \leq d}$ (the discounted process). \\
- some financial assets (as options). Each of these financial assets is modeled  by an essentially bounded  ${\cal F}_{\tau_l}$-measurable function  $Y^l$ ($\tau_l$ the  maturity date can be a stopping time). We assume furthermore that  at time zero,   a bid price and  ask price are observed in the market for any of the assets $Y^l$.\\   Define now the notion of admissibility with respect to the reference family $((S^k)_{0 \leq k \leq d},(Y^l)_{1 \leq l \leq p})$. In \cite{BN02} we have introduced the notion of admissibility for a static pricing function. We extend here this notion to a dynamic context.\\  
   For any  $l$ denote ${C^l}_{ask}$ ( resp ${C^l}_{bid}$) the observed ask (resp bid) price observed in  the market at time $0$ for the asset $Y^l$.
 We give the following definition of strong admissibility:

\begin{definition}
A Dynamic  Pricing Procedure  $(\Pi_{\sigma, \tau})_{0 \leq \sigma \leq \tau}$ is strongly admissible with respect to the 
reference family $((S^k)_{0 \leq k \leq d},(Y^l)_{1 \leq l \leq p})$  and the observed bid and  ask prices $({C^l}_{bid},{C^l}_{ask})_{1 \leq l \leq p}$ if   \begin{itemize}
\item it extends the dynamics of the process $(S^k)_{0 \leq k \leq d}$.
\item it is compatible with the observed bid and ask prices for the $(Y^l)_{1 \leq l \leq p}$:
\begin{equation}
\forall \; 1 \leq l \leq p \;\; \;{C^l}_{bid} \leq -\Pi_{0,\tau_l}(-Y^l) \leq\; 
\Pi_{0,\tau_l}(Y^l) \leq {C^l}_{ask}\;\;
\label{A2}
\end{equation}
\end{itemize}
\label{definition2}
\end{definition}
 The following proposition gives a characterization of  strongly admissible TCPP. 
\begin{proposition}
A TCPP   is strongly admissible with respect to the 
reference family $((S^k)_{0 \leq k \leq d},(Y^l)_{1 \leq l \leq p})$  and the observed bid and ask prices\\ $({C^l}_{bid},{C^l}_{ask})_{1 \leq l \leq p}$ if and only if 
\begin{itemize}
\item Any probability measure $R \ll P$ such that $\alpha^m_{0,\infty}(R) < \infty$ is a local martingale measure with respect to any process  $S^k$.
\item For any probability measure $R \ll P$, for any stopping time $\tau$,
\begin{equation} 
\alpha^m_{0,\tau}(R) \geq \sup(0,\sup_{\{l \;|\; \tau_l \leq \tau\}}(C^l_{bid}-E_R(Y^l),E_R(Y^l)-C^l_{ask})
\label{eqSA}
\end{equation}
\end{itemize}
\label{propSA}
\end{proposition}
{\bf Proof}:  This follows easily from the characterization of TCPP extending the dynamics of $S^k$ (Proposition \ref{propEMM}) and from the dual representation of $\Pi_{0,\tau}$, equation (\ref{eq_4}).

\begin{remark}
From Theorem \ref{thmarbitrage}, and Proposition \ref{propSA}, it follows that the existence of a strongly admissible No Free Lunch TCPP relies on the existence of an equivalent local martingale measure $Q_0$ for the process $(S^k)_{0 \leq k \leq d}$, such that 
\begin {equation}
\forall \; 1 \leq l \leq p \;\; \;{C^l}_{bid} \leq E_{Q_0}(Y^l) \leq {C^l}_{ask} 
\label{eqSAEMM}
\end{equation}
\end{remark}
A detailled study of TCPP calibrated with observed bid ask prices for options is the subject of \cite{BN05} (in preparation). In particular it is shown in \cite{BN05} that the existence of  an equivalent local martingale measure $Q_0$ satisfying the inequalities (\ref{eqSAEMM})
is equivalent to the No Free Lunch condition for the reference family 
 $((S^k)_{0 \leq k \leq d},(Y^l)_{1 \leq l \leq d},( C^l_{bid}, C^l_{ask})_{1 \leq l \leq d}$. (This result  generalizes the usual Kreps Yan Theorem).\\
Strong admissible No Free Lunch TCPP lead to sharper bounds.
Denote ${\cal M}^e$ the set of equivalent local martingale measures for the process $(S^k)_{0 \leq k \leq d}$. For any $Q\in {\cal M}^e$ denote 
$\beta(Q)= \sup(0,\sup_{1 \leq l \leq p}(C^l_{bid}-E_Q(Y^l),E_Q(Y^l)-C^l_{ask}))$.
\begin{proposition}
 For any No Free Lunch TCPP strongly admissible with respect to the reference family $((S^k)_{0 \leq k \leq d},(Y^l)_{1 \leq l \leq p})$  and the observed bid and ask prices $({C^l}_{bid},{C^l}_{ask})_{1 \leq l \leq p}$, for any 
financial product $X \in L^{\infty}({\cal F}_{\infty})$, 
\begin{equation}
\inf_{Q \in {\cal M}^e} (E_Q(X)+\beta(Q))\leq -\Pi_{0,\infty}(-X) \leq \Pi_{0,\infty}(X) \leq \sup_{Q  \in {\cal M}^e} (E_Q(X)-\beta(Q))
\label{eqTB}
\end{equation}
\end{proposition}
 {\bf Proof}: 
 This result follows from the dual representation for No Free Lunch TCPP, Corollary \ref{corNFL} and the Proposition \ref{propSA} (applied with $\tau=\infty$). \\

 Notice that the bounds in inequation (\ref{eqTB}) are tighter than the bounds obtained from surreplication (and subreplication), as soon as the family of options taken into account is sufficiently rich. Indeed as soon as one option is one of the options of the reference family, the bonds of any TCPP are inside the interval $[C^l_{bid},{C^l}_{ask}]$ corresponding to the observed bid and ask prices in the market. Furthermore adding some new option to the reference family is equivalent to have  more conditions in condition (\ref{eqSA}) and therefore leads to no wider and sometimes  tighter bounds. We can notice that these new tighter bounds are just obtained from the information already contained in the market. There is no need to introduce some external restriction or preference. In  \cite{BN05} we provide examples of strongly admissible No Free Lunch TCPP, making use of the theory of right continuous  BMO martingales. 
 
\section{Conclusion}
 We have introduced an axiomatic approach for bid ask dynamic pricing in general financial markets.
 We have defined the notion of TCPP (time consistent dynamic pricing procedure continuous from below) which assigns to any financial product both an ask  price process and a bid price process in such a way that it takes into account both transaction costs and liquidity risk. Therefore the ask price process $\Pi_{\sigma,\tau}(X)$ is a convex function of $X$.\\
 We have defined two  notions of No Free Lunch for a pricing procedure, a static one and a dynamic one. For any time consistent dynamic pricing procedure the two notions are equivalent. Our main result (generalizing that of Jouini and Kallal \cite{JK})  is that the No Free Lunch property for a TCPP is equivalent to the existence of an equivalent probability measure $R$ such that for any $X \in L^{\infty}(\Omega,{\cal F}_{\infty},P)$, the martingale $E_R(X | {\cal F}_{\sigma})$ is always between the Bid Price Process $-\Pi_{\sigma,\tau}(-X)$ and the Ask Price Process $\Pi_{\sigma,\tau}(X)$. Furthermore the ask price process associated with any financial instrument is a $R$-supermartingale process which has a cadlag modification (when the No Free Lunch property is satisfied). \\
  In the last section \ref{sec:admissible} we considered   the usual setting for financial markets, where the dynamics of some reference assets is assumed to be known $(S^k)_{0 \leq k \leq p}$. We have proved that  No Free Lunch TCPP extending the dynamics  of the reference assets have a dual representation in terms of  equivalent local  martingale measures for the reference process $(S^k)_{0 \leq k \leq p}$.
 TCPP are usefull in order to produce bid ask spreads for options close to the observed bid ask spreads in real financial markets. A family of TCPP well adapted to financial markets, is the TCPP calibrated with observed bid and ask prices for options. It allows for the   production for  market makers of bid and ask prices for any  new financial instrument. Furthermore the prices produced that way are consistent with the bid ask spreads that one wants to take into account and lead for other products to tighter bounds than the usual ones.  
 The systematic study of TCPP extending the dynamics of some reference assets and compatible with observed bid ask prices for some options is the subject of a forthcoming paper \cite{BN05}.

\end{document}